\documentclass[11pt]{article} 

\usepackage[utf8]{inputenc} 

\usepackage{amsmath}
\usepackage{amssymb}
\usepackage{amsthm}
\usepackage{tikz}
\usetikzlibrary{arrows}

\newtheorem*{thm1}{Theorem 1}
\newtheorem*{cor2}{Corollary 2}
\newtheorem*{thm3}{Theorem 3}
\newtheorem*{thm4}{Theorem 4}

\usepackage{geometry} 
\usepackage{graphicx} 

\usepackage{booktabs} 
\usepackage{array} 
\usepackage{paralist} 
\usepackage{verbatim} 
\usepackage{subfig} 

\usepackage{fancyhdr} 
\usepackage{sectsty}
\allsectionsfont{\sffamily\mdseries\upshape} 

\usepackage[nottoc,notlof,notlot]{tocbibind} 
\usepackage[titles,subfigure]{tocloft} 

\title{Local Maxima of Quadratic Boolean Functions}
\author{Hunter Spink}

\begin{document}
\maketitle
\begin{abstract}
How many strict local maxima can a real quadratic function on $\{0,1\}^n$ have? Holzman conjectured a maximum of $n \choose \lfloor n/2 \rfloor$. The aim of this paper is to prove this conjecture. Our approach is via a generalization of Sperner's theorem that may be of independent interest.
\end{abstract}
\section{Introduction}
Let $\Theta$ be a real quadratic function (polynomial of total degree $\le 2$) in $n$ variables $x_1,\ldots x_n$. A \textit{strict local maximum} (or just \textit{local maximum}) of $\Theta$ on the discrete cube $Q_n=\{0,1\}^n$ is a point whose value is strictly larger than all of its neighbors.
As we are only concerned with the value of $\Theta$ on $Q_n$, we may assume when convenient that all terms are of degree 2, as the constant term is irrelevant, and we can replace $x_i$ with $x_i^2$ if necessary.
In this paper, we prove the following conjecture attributed to Ron Holzman (see \cite[p.~3]{problem}):
\begin{thm1}
Let $\Theta$ be a quadratic function on $Q_n$. Then $\Theta$ has at most $n \choose \lfloor n/2 \rfloor$ local maxima.
\end{thm1}
This bound is attained for example when $\Theta=-(x_1+\ldots +x_n-\lfloor n/2 \rfloor)^2$.

From this, we can deduce the following:
\begin{cor2}
A (possibly degenerate) parallelepiped in $\mathbb{R}^n$ can have at most $n \choose \lfloor n/2 \rfloor$ vertices that are strictly closer to the origin than all of their neighbors.
\end{cor2}
\begin{proof}
Let $\Theta$ be the form $-\sum x_i^2$; we are counting strict local maxima of $\Theta$ on the parallelepiped. As being a strict local maximum is clearly an open condition, we may perturb the sides so that the parallelepiped is not degenerate. There is an affine transformation $\tau$ taking this parallelepiped to $Q_n$, so composing $\Theta$ with $\tau^{-1}$, we are done by Theorem 1.
\end{proof}
[We remark in passing that Corollary 2 implies the Littlewood-Offord theorem of Kleitman \cite{kleitman}. Indeed, the theorem is equivalent to showing that the number of vertices of a parallelepiped with all side lengths at least $2$ that can land in the interior of a disc of radius $1$ is at most $n \choose n/2$. The result follows by noting that a vertex landing inside the disc must have all neighbors outside of the disc, hence farther from the disc's center.]

Let $[n]=\{1,2,\ldots n\}$, and let $\Delta$ denote symmetric set difference. We define the upper $i$th level of a family of subsets $\mathcal{F}$ of $[n]$ as $$\mathcal{F}^{i+}=\{A \subseteq [n] : i \not\in A, A \cup \{i\} \in \mathcal{F}\},$$ and the lower $i$th level as $$\mathcal{F}^{i-}=\{A \subseteq [n] : i \not\in A, A \in \mathcal{F}\}.$$
Finally, for a subset $S \subseteq [n]$, we define $\mathcal{F}\Delta S$ be the family of sets $A \Delta S$ for all $A \in \mathcal{F}$.

We will find that a key step is to prove the following combinatorial result:
\begin{thm3}
Given $S_i \subseteq [n]$ for $i=1,2,\ldots, n$,  and a family of subsets $\mathcal{F}$ such that, for each $i$, no element of $(\mathcal{F}\Delta S_i)^{i+}$ contains an element of $(\mathcal{F} \Delta S_i)^{i-}$, then $|\mathcal{F}| \le {n \choose\lfloor n/2\rfloor}$.
\end{thm3}
Theorem 3 is perhaps interesting in its own right, because it gives a generalization of Sperner's theorem (which is the case when all $S_i=\emptyset$).

At the end of the paper, we analyze the cases of when the form attains the maximum number of local maxima. It turns out our method allows us not only to deduce the structure of the quadratic function when we attain equality, but also when we are within $\frac{1}{n}{n \choose \lfloor n/2 \rfloor}$ of the optimal solution $n \choose \lfloor n/2 \rfloor$. We do not know how close this bound of $\frac{1}{n}{n \choose \lfloor n/2 \rfloor}$ is to being optimal.

This paper is self-contained. See \cite{bela} for general background on set systems and Sperner's theorem.
\section{Proof}
 We frequently view the discrete cube as the family of all subsets of $[n]$ in the obvious way (using the $x_i$ as indicator functions). Assuming all terms of $\Theta$ are of degree 2 (so that it is a quadratic form), we have its associated symmetric matrix $(q_{ij})$, where $\Theta(x_1,\ldots, x_n)=\sum_{i,j} q_{ij}x_ix_j$.

There are $2^n$ graph automorphisms one gets by taking some subset $S$ of $[n]$, and replacing each $A \subseteq [n]$ by $A\Delta S$: we preserve the quadratic form by replacing $x_i$ by $1-x_i$ for all $i \in S$. When we apply such an automorphism we will say we are ``changing the origin", or ``changing the base", since this is how one can view such an action geometrically (after suitable changes of signs).

\subsection{Proof that Theorem 3 implies Theorem 1}
By perturbing the form slightly, we may assume all $q_{ij}$ are nonzero.
We note first that changing the origin to $B=(b_1,\ldots, b_n)$ has the effect on the off-diagonal entries of $(q_{ij})$ of flipping the signs in the $i$th row and column for all $i$ with $b_i=1$ (so leaving unchanged any off-diagonal $q_{ij}$ with $b_i=b_j=1$). The behaviour of the diagonal entries is more complicated, and does not need to be considered.

The difference between the value of $\Theta$ on the $x_i=1$ plane and the $x_i=0$ plane as a function of the remaining coordinates is given by the linear function $2\sum_{j \ne i} q_{ij}x_i+q_{ii}$. Changing the origin to $S_i$, where $S_i$ is the set of all $j \ne i$ with $q_{ij}$ positive, we can assume $q_{ij}<0$ for all $j \ne i$. In this new coordinate system associated to $S_i$, we cannot have $B$ in the $x_i=1$ plane containing $C$ in the $x_i=0$ plane, with both $B,C$ local maxima. This is because $B$ being a local maximum means we must have $2\sum_{j \ne i} q_{ij}b_j+q_{ii}>0$ and $C$ being a strict local maximum means $2\sum_{j \ne i} q_{ij}c_j+q_{ii}<0$, so we get $2\sum_{j \ne i} q_{ij} (b_j-c_j)>0$, which is clearly false.

Also, $B$ containing $C$ in the $S_i$-coordinate system is equivalent in the original coordinate system to the statement $B \Delta S_i \supseteq C \Delta S_i$. Hence we have reduced to Theorem 3.
\subsection{Proof of Theorem 3}
Inspired by Kleitman's proof of the Littlewood-Offord problem \cite{kleitman} (or see \cite[Ch.4]{bela} for general background), we seek a ``symmetric quasichain decomposition" of $[n]$ (described below), where a ``quasichain" will be a family with some property $P$ which implies at most one member of $\mathcal{F}$ lies inside it. The heart of this proof is the definition of a ``quasichain", which allows this method to go through. The definition is rather surprising and seems contrived, however as we will see, once we have this definition the proof is straightforward.

A ``symmetric quasichain decomposition" in this case will then be a decomposition of the $n$-cube inductively built up from the $1$-cube by taking the quasichain decomposition of the $k$-cube, duplicating each quasichain with $k+1$ added to each set, then removing precisely one set from each duplicate and adding it to the original in such a way that the new families formed remain quasichains. It is not hard to prove (or see \cite[p.~17-20]{bela}) that this process will result in exactly $n \choose \lfloor n/2 \rfloor$ quasichains, as required.

For sets $B,C \subseteq [n]$, we write $B \supseteq_{S_i} C$ to mean $B \Delta S_i \supseteq C \Delta S_i$.

We define a \textit{quasichain} to be a colored tournament (with colors in $[n]$) on a family of subsets $\mathcal{G}=\{G_1,\ldots G_k\}$ of $[n]$, such that

\begin{enumerate}[i)]
\item Whenever there is a directed edge from $G_x$ to $G_y$ of color $i$, then $i \in G_x\Delta S_i$, $i \not \in G_y\Delta S_i$, and $ G_x \supseteq_{S_i} G_y$.
\item For any subset of the colors (including the empty set), if we swap the direction of edges associated to those colors, then the resulting tournament is acyclic (or equivalently transitive).
\end{enumerate}

It is easy to check that the acyclicity condition is equivalent to saying that no triangle has 3 distinct colors, any monochromatic triangle is acyclic, and any triangle with 2 distinct colors has the 2 edges with the same color either both leaving, or both entering, the same vertex.

Note that a quasichain does \textit{not} contain all possible information about $\supseteq_{S_i}$ contaiment between its various members, but rather remembers only one such containment for every pair (just enough information to guarantee that at most one element from the pair can be a local maximum from the condition in the theorem).

We write $G_x \to^i G_y$ when there is an edge from $G_x$ to $G_y$ of color $i$. Sometimes we will reduce to the case $i \not \in S_i$, in which case $G_x \to^i G_y$ implies $i \in G_x$, $i \not\in G_y$, and $G_x \supseteq_{S_i} G_y$.
\begin{proof}
Note that for a fixed family $\mathcal{F}$, the hypothesis of the theorem is easily checked to be invariant under base-change by an arbitrary subset $A$ (i.e. taking $F \to \mathcal{F}\Delta A$ and $S_i \to S_i \Delta A$). It is also invariant replacing $S_i$ by $S_i^c$ for any $i$, so in particular we may assume that $i \not\in S_i$ for all $i$. Also note that given a quasichain $\mathcal{G}$, it remains a quasichain after replacing $S_i$ by $S_i^c$ if one swaps the directions of all $i$-colored edges (this is the reason we need the acyclicity condition in the induction hypothesis, so that complementing the $S_i$ preserves the property of being a quasichain).

If we have a quasichain $\mathcal{G}$, then base-changing to $A$ followed by complementing every $S_i$ for which  $i \in S_i$ has the following effect: each $S_i$ is changed to $S_i\Delta A$ if $i \not\in A$ or $(S_i \Delta A)^c$ otherwise, $G_i$ turns into $G_i\Delta S_i$, and the direction of all $i$-colored arrows are swapped if $i \in A$. (We remark for motivation that this net effect is equivalent to what happens when we base-change by $A$ in the quadratic form case, where $S_j=\{i \ne j : q_{ij}>0\}$, as the complementation is ``built in" to the definition of $S_j$).

Given $S_1,\ldots, S_n$, assume by induction we have a symmetric quasichain decomposition for $Q_{n-1}$ (using colors in $[n-1]$) associated to the sets $S_i \cap \{1,\ldots, n-1\}$ for $i=1,\ldots, n-1$.  We will produce from this a symmetric quasichain decomposition for $Q_n$ (using colors in $[n]$) associated to the sets $S_i$ for $i=1,\ldots, n$. By base-changing with resect to $\{i : n \in S_i\}$, and complementing as described above, we can assume without loss of generality that $i \not \in S_i$ and $n \not\in S_i$, for $i=1,\ldots,n$.

Note that if we treat a quasichain $\mathcal{C}$ in $Q_{n-1}$  as a quasichain in $Q_n$, then it remains a quasichain (obviously), and if we add $n$ to each set in $\mathcal{C}$ (denoted $\mathcal{C}+n$), then this is also a quasichain. By the standard symmetric quasichain construction discussed previously, it suffices to show that there exists an element of $\mathcal{C}+n$ which we can remove from it and add to $\mathcal{C}$ such that the two newly constructed directed graphs are quasichains.

As $\mathcal{C}$ is acyclic, it has a maximal element $A$. A subgraph of a quasichain is clearly a quasichain, so $\mathcal{C}+n-\{A+n\}$ is a quasichain (here minus means we remove $A+n$ from the quasichain). Thus it suffices to show that $\mathcal{C}+\{A+n\}$ is in fact a quasichain once we appropriately color edges containing $A+n$.

If $A \supseteq_{S_i} B$, then since $n \not\in S_i$, we have $A+n \supseteq_{S_i} A \supseteq_{S_i} B$, so by transitivity of $\supseteq_{S_i}$, we have $A+n \supseteq_{S_i} B$. Thus if $A \to^i B$, then $i\not\in B$, and $i \in A$, so clearly $i \in A+n$, and we may set $A+n \to^i B$.
Also, as both $n \not\in S_n$ and $n \not\in A$, we may also set $A+n \to^n A$.

Thus it suffices to show that the newly constructed graph satisfies the acyclicity condition. After swapping some directions associated to a subset of the colors, we have a tournament $H$, with vertices $x$, $y$ corresponding to $A$, $A+n$ respectively. For all $z$ then, we have either $x \to z$, $y\to z$, or $z\to x$, $z \to y$. If we have a cycle, then identifying $x$ and $y$ yields a cycle in the original graph (since $x$ and $y$ have the same incoming/outgoing edges, this identification is well-defined), which is a contradiction.
\end{proof}
\tikzset
{
	e1/.style=
	{
		red
	},
	e2/.style=
	{
		blue
	},
	e3/.style=
	{
		black
	},
	e4/.style=
	{
		purple
	},
	e5/.style=
	{
		orange
	},
}
\section{Analysis of the quadratic form close to equality}

From now on, we are working with the $S_j$ associated to a quadratic form $\Theta$ (recall we defined $S_j=\{i \ne j : q_{ij}>0\}$). Then $i \not \in S_i$, $i \in S_j$ if and only if $j \in S_i$, and these properties are invariant under base-change. Base-change by $S_1$ so that $S_1=\emptyset$. We will show that if all $S_i$ are not the empty set (i.e. in this coordinate system there is an off-diagonal entry which is positive), then we are bounded away from the optimal solution by a factor of $\frac{1}{n}$. Indeed, suppose without loss of generality $S_2$ contains 3. Then $S_3$ contains 2, so the intersections with $\{1,2,3\}$ of $S_1,S_2,S_3$ are $\emptyset, \{3\},\{2\}$ respectively. The first three stages yield the quasichain decomposition below (in bold lines):

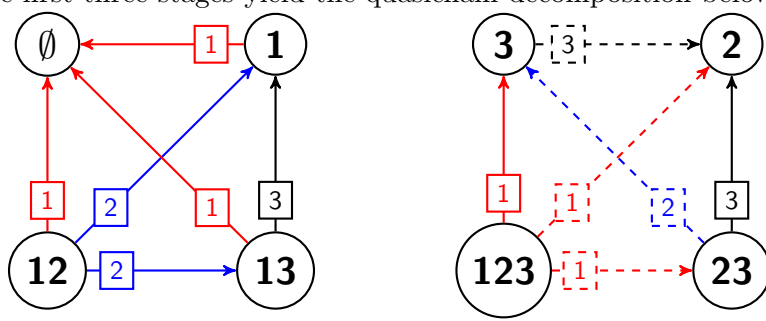
\begin{tikzpicture}[->,>=stealth',node distance=3cm,
  thick,main node/.style={circle,draw,font=\sffamily\Large\bfseries}]

  \node[main node] (0) {$\emptyset$};
  \node[main node] (1) [right of=0] {1};
  \node[main node] (12) [below of=0] {12};
  \node[main node] (13) [right of=12] {13};
  \node[main node] (3) [right of=1] {3};
  \node[main node] (123) [below of=3] {123};
  \node[main node] (2) [right of=3] {2};
  \node[main node] (23) [below of=2] {23};

  \path[every node/.style={rectangle,draw, fill=white, font=\sffamily\small, pos=0.2}]
	(1)	edge[e1] node {1} (0)
	(12)	edge[e1] node {1} (0)
		edge[e2] node {2} (1)
		edge[e2] node {2} (13)
	(13)	edge[e3] node {3} (1)
		edge[e1] node {1} (0);
    \path[every node/.style={rectangle,draw, fill=white, font=\sffamily\small, pos=0.2}]
	(3)	edge[e3,dashed] node {3} (2)
	(123)	edge[e1] node {1} (3)
		edge[e1,dashed] node {1} (2)
		edge[e1,dashed] node {1} (23)
	(23)	edge[e2,dashed] node {2} (3)
		edge[e3] node {3} (2);
\end{tikzpicture}

Using the dotted lines, we see that we can ``glue" together two of these quasichains to make a single quasichain. We want to understand the ``evolution" of this quasichain as it goes through the symmetric chain algorithm. We know that a 2-quasichain after $k$ steps becomes $k+1 \choose \lfloor (k+1)/2 \rfloor$ quasichains. As it evolves into a 4-quasichain plus two 2-quasichains after two steps, if $g(k)$ is the number of quasichains the 4-quasichain evolves into after $k$ steps, it satisfies
$${k+3 \choose \lfloor (k+3)/2 \rfloor}=2{k+1 \choose \lfloor (k+1)/2 \rfloor}+g(k).$$
After we reach $n$-dimensions therefore, these two 4-quasichains will have evolved to
$$2g(n-3)=2{n \choose \lfloor n/2 \rfloor}-4{n-2 \choose \lfloor (n-2)/2 \rfloor}.$$
The difference between the actual bound and this is
$$4{n-2 \choose \lfloor (n-2)/2 \rfloor}-{n \choose \lfloor n/2 \rfloor},$$
which is equal to $$\frac{1}{n-c}{n \choose \lfloor n/2 \rfloor}$$
where $c$ is 0 or 1 depending on whether $n$ is odd or even respectively.

If every $S_i$ is in fact empty, then the condition the $S_i$ create is precisely the normal antichain condition, so the family creates an actual antichain.
Thus we have the following:
\begin{thm4}
If the number of local maxima of a quadratic function on $Q_n$ is greater than $(1-\frac{1}{n-c}){n\choose n/2}$, where $c=0,1$ if $n$ is odd/even respectively, then there is an automorphism of $Q_n$ such that after applying it, the local maxima form an antichain.
\end{thm4}

\section{Acknowledgements}
I would like to thank Imre Leader for helpful conversations. I would also like to thank Trinity College, Cambridge for providing me with financial support.

\end{document}